\documentclass[12pt]{amsart}
\usepackage{amssymb,amsmath,amsfonts,latexsym,setspace}
\usepackage{bm}
\newcommand{\newsection}[1]{\setcounter{equation}{0} \section{#1}}
\newcommand{\bea}{\begin{eqnarray}}
\newcommand{\eea}{\end{eqnarray}}

\newcommand{\vp}{\varphi}

\newcommand{\clb}{\mathcal{B}}

\newcommand{\cle}{\mathcal{E}}

\newcommand{\clg}{\mathcal{G}}
\newcommand{\clh}{\mathcal{H}}
\newcommand{\clk}{\mathcal{K}}

\newcommand{\clm}{\mathcal{M}}

\newcommand{\clo}{\mathcal{O}}

\newcommand{\clq}{\mathcal{Q}}

\newcommand{\cls}{\mathcal{S}}

\newcommand{\raro}{\rightarrow}

\newcommand{\G}{\Gamma}
\newcommand{\p}{\partial}

\def \qed {\hfill \vrule height6pt width 6pt depth 0pt}
\def\textmatrix#1&#2\\#3&#4\\{\bigl({#1 \atop #3}\ {#2 \atop #4}\bigr)}
\def\dispmatrix#1&#2\\#3&#4\\{\left({#1 \atop #3}\ {#2 \atop #4}\right)}
\newcommand{\be}{\begin{equation}}
\newcommand{\ee}{\end{equation}}
\newcommand{\ben}{\begin{eqnarray*}}
\newcommand{\een}{\end{eqnarray*}}

\newcommand{\NI}{\noindent}

\newcommand{\bi}{\begin{itemize}}
\newcommand{\ei}{\end{itemize}}

\newtheorem{Theorem}{\sc Theorem}[section]
\newtheorem{Lemma}[Theorem]{\sc Lemma}
\newtheorem{Proposition}[Theorem]{\sc Proposition}
\newtheorem{Corollary}[Theorem]{\sc Corollary}
\newtheorem{Definition}[Theorem]{\sc Definition}
\newtheorem{Example}[Theorem]{\sc Example}
\newtheorem{Remark}[Theorem]{\sc Remark}
\newtheorem{Note}[Theorem]{\sc Note}
\newtheorem{Question}{\sc Question}
\newtheorem{ass}[Theorem]{\sc Assumption}
\newcommand{\bt}{\begin{Theorem}}
\def\beginlem{\begin{Lemma}}
\def\beginprop{\begin{Proposition}}
\def\begincor{\begin{Corollary}}
\def\begindef{\begin{Definition}}
\def\beginexamp{\begin{Example}}
\def\beginrem{\begin{Remark}}
\def\beginq{\begin{Question}}
\def\beginass{\begin{ass}}
\def\beginnote{\begin{Note}}
\newcommand{\et}{\end{Theorem}}
\def\endlem{\end{Lemma}}
\def\endprop{\end{Proposition}}
\def\endcor{\end{Corollary}}
\def\enddef{\end{Definition}}
\def\endexamp{\end{Example}}
\def\endrem{\end{Remark}}
\def\endq{\end{Question}}
\def\endass{\end{ass}}
\def\endnote{\end{Note}}

\begin{document}

\title[Similarity of Quotient Hilbert modules in the Cowen-Douglas Class]{Similarity of Quotient Hilbert modules in the Cowen-Douglas Class}

\author[Kui Ji]{Kui Ji}
\address{Department of Mathematics, Hebei Normal University, Shijiazhuang, Hebei 050016, China}
\email{jikui@hebtu.edu.cn}

\author[Jaydeb Sarkar]{Jaydeb Sarkar}
\address{Indian Statistical Institute, Statistics and Mathematics Unit, 8th Mile, Mysore Road, Bangalore, 560059, India}
\email{jay@isibang.ac.in, jaydeb@gmail.com}




\begin{abstract}
In this paper, we consider the similarity and quasi-affinity
problems for Hilbert modules in the Cowen-Douglas class associated
with the complex geometric objects, the hermitian anti-holomorphic
vector bundles and curvatures. Given a ``simple'' rank one
Cowen-Douglas Hilbert module $\clm$, we find necessary and
sufficient conditions for a class of Cowen-Douglas Hilbert modules
satisfying some positivity conditions to be similar to $\clm \otimes
\mathbb{C}^m$. We also show that under certain uniform bound
condition on the anti-holomorphic frame, a Cowen-Douglas Hilbert
module is quasi-affinity to a submodule of the free module $\clm
\otimes \mathbb{C}^m$.
\end{abstract}

\maketitle

\newsection{Introduction}
One of the most challenging problems in operator theory is to
determine when two given bounded linear operators are similar. More
precisely, let $T$ and $R$ be two bounded linear operators on
Hilbert spaces $\clh$ and $\clk$, respectively. When does there
exists an invertible bounded linear map $X : \clh \raro \clk$ such
that $X T = R X$?

There are many fascinating subtleties connected with the similarity
problem (see \cite{JW}, \cite{P}, \cite{Stam}, \cite{NF},
\cite{NF2}). However, the problem becomes more tractable if one
impose additional assumptions on the operators. In particular, there
are several characterizations of operators similar to unitaries or
isometries or even contractions.

In \cite{U}, Uchiyama proposed a characterization for contractions
in the Cowen-Douglas class which are similar to the adjoint of the
multiplication operators on the Hardy space with finite
multiplicity. One of the main tools used in the work by Uchiyama is
the structure of the tensor product bundle corresponding to a given
hermitian holomorphic vector bundle. Later, Kwon and Treil \cite{KT}
found some additional characterizations which involves the
curvature, in the sense of Cowen-Douglas, and the Carleson measure
\cite{C} of the underlying operators (see also Douglas, Kwon and
Treil \cite{DKT} in the setting of $n$-hypercontractions \cite{Ag}).

In the present study, we set up the similarity problem in a more
general framework and provide some characterizations of operators in
the Cowen-Douglas class which are similar to the adjoint of the
multiplication operators on standard reproducing kernel Hilbert
spaces of holomorphic functions (like weighted Bergman spaces). We
prove that the earlier characterizations of operators similar to the
adjoint of multiplication operators are valid beyond the class of
contractions and $n$-hypercontractions. In particular, our results
includes the similarity problem for the weighted Bergman spaces with
not necessarily integer weights. Our framework is based on the
$\frac{1}{K}$-calculus introduces by Arazy and Englis \cite{AE}. The
similarity results of this paper are significantly more general than
those obtained in \cite{DKT}, \cite{KT} and \cite{U}. Moreover, the
study of weighted Bergman spaces with not necessarily integer
weights appears to be more fruitful and rewarding from analytic,
geometric and representation theoretic point of views (for instance,
see \cite{KM} and \cite{MU}).

We now summarize the content of this paper. In Section 2, we give a
self-contained presentation of the theory of Cowen-Douglas Hilbert
modules in the language of reproducing kernel Hilbert modules. The
next section is devoted to assembling materials, like a projection
formula, derivatives of holomorphic maps and curvature, which we
will use in the sequel, from various sources. Here, however, our
approach will be completely new. In Section 4, we develop the notion
of Cowen-Douglas atoms. In Section 5, we obtain results concerning
tensor product bundles and quotient modules. Section 6 deals with
the quasi-affinity properties of Cowen-Douglas Hilbert modules. In
Section 7, we relate the curvatures to the derivatives of
holomorphic maps and obtain some similarity classification results.
The final section lists some open problems.

\newsection{Cowen-Douglas Hilbert Modules}

In this section we introduce the basic concepts and known results
related to the Cowen-Douglas class \cite{CD}.

Let $T$ be a bounded linear operator on a Hilbert space $\clh$. Then
$\clh$ is a module over $\mathbb{C}[z]$ in the following sense:
\[
p \cdot f \mapsto p(T) f \quad \quad \quad \quad (p \in
\mathbb{C}[z], f \in \clh).
\]
The above module is usually called the \textit{Hilbert module} over
$\mathbb{C}[z]$ (see \cite{DP}).

Note that a Hilbert module $\clh$ over $\mathbb{C}[z]$ is uniquely
determined by the underlying operator $T$ via the module
multiplication operator by the coordinate function $z$:
\[
M_z f : = z \cdot f = T f, \quad \quad \quad (f \in \clh),
\]
and vice versa. We say that $\clh$ is \textit{contractive} if
\[
I_{\clh} - M_z M_z^* \geq 0.
\]

We denote the space of all bounded linear operators from a Hilbert
space $\clh$ to another Hilbert space $\clk$ by $\clb(\clh, \clk)$,
and by $\clb(\clh)$ if $\clk = \clh$.

Let $\clh_1$ and $\clh_2$ be two Hilbert modules over
$\mathbb{C}[z]$. Then $M \in \clb(\clh_1, \clh_2)$ is said to be a
module map if $M(p \cdot f)=p \cdot (M f)$ for all $p  \in
\mathbb{C}[z]$ and $f \in \clh_1$.

Now we recall the definition of the Cowen-Douglas class \cite{CD}.

\begin{Definition}
Let $m$ be a positive integer, and let $\Omega$ be a domain in
$\mathbb{C}$. The Cowen-Douglas class on $\Omega$ of rank $m$,
denoted by $B_m(\Omega)$ is the set of all Hilbert modules $\clh$
over $\mathbb{C}[z]$ such that

(i) $\sigma(M_z^*) \subseteq \Omega$,

(ii) $\mbox{ran} (M_z - wI_{\clh})^* = \clh$ for all $w \in \Omega$,

(iii) $\mbox{dim~}\mbox{ker}(M_z - wI_{\clh})^* = m$ for all $w \in
\Omega$, and

(iv) $\overline{\mbox{span}}\{ \mbox{ker} (M_z - w I_{\clh})^* : w
\in \Omega\} = \clh$.
\end{Definition}

A Hilbert module $\clh$ is said to be a \textit{Cowen-Douglas
Hilbert module} if $\clh \in B_m(\Omega)$ for some positive integer
$m$.


A Cowen-Douglas Hilbert module $\clh$ in $B_m(\Omega)$ defines a
hermitian anti-holomorphic vector bundle $E_{\clh}$ over $\Omega$
where
\[E_{\clh} = \{(w, h) \in \Omega \times \clh : M_z^* h = \bar{w}
h\},\]with the projection map $\pi_{\clh} : E_{\clh} \raro \Omega$
defined by $\pi_{\clh}(w, h) = h$ for all $w \in \Omega$ and $h \in
\clh$. More precisely, $E_{\clh}$ is the anti-holomorphic vector
bundle implemented by the anti-holomorphic map $w \mapsto
E_{\clh}(w) := \mbox{ker~} (M_z - w I_{\clh})^*$, the pull-back
bundle of the Grassmannian $GF(m, \clh)$ (see \cite{CD}) and hence
locally at each point $w \in \Omega$, there exists anti-holomorphic
$\clh$-valued functions $\{\gamma_{i,w} : 1 \leq i \leq m\}$ such
that
\[
\mbox{span} \{\gamma_{i,w} : 1 \leq i \leq m\} = \mbox{ker} (M_z -
w)^*.
\]
Also it follows from a theorem of Grauert \cite{grauert} that
$\gamma_{i,w}$ can be defined on all of $\Omega$.

The rigidity theorem (Theorem 2.2 in \cite{CD}) states that a pair
of Hilbert modules $\clh$ and $\tilde{\clh}$ in $B_m(\Omega)$ are
unitarily equivalent if and only if the corresponding hermitian
anti-holomorphic vector bundles $E_{\clh}$ and $E_{\tilde{\clh}}$
are equivalent.

We now briefly recall the notion of a reproducing kernel Hilbert
module, which in turn is closely related to the Cowen-Douglas
Hilbert modules. Let $\cle$ be a Hilbert space. A Hilbert module
$\clh \subseteq \mathcal{O}(\Omega, \cle)$, where
$\mathcal{O}(\Omega, \cle)$ is the space of $\cle$-valued
holomorphic functions on $\Omega$, is said to be a
\textit{reproducing kernel Hilbert module} if

(i) the evaluation map $ev_w : \clh \raro \cle$, $w \in \Omega$,
defined by $ev_w(f) = f(w)$, $f \in \clh$, is bounded, and

(ii) the module multiplication operator $M_z$ is given by the
multiplication operator by the coordinate function $z$.

Let $\clh \in B_m(\Omega)$ with an anti-holomorphic frame
$\{\gamma_{i, w} : 1 \leq i \leq m\}$ of $E_{\clh}$, and let $G(z,
w)$ be the corresponding Gram matrix
\[
G(z, w) = (\langle \gamma_{j, w}, \gamma_{i, z} \rangle_{\clh} )_{m
\times m}, \quad \quad (z, w \in \Omega).
\]
Define $K : \Omega \times \Omega \raro \clm_m(\mathbb{C})$ by
\[
K(z, w) = G(z, w) \quad \quad (z, w \in \Omega).
\]
Then $K$ is a positive definite kernel, and the corresponding
reproducing kernel Hilbert space $\clh_K \subseteq \clo(\Omega,
\mathbb{C}^m)$ is a reproducing kernel Hilbert  module. Define $X :
\clh \raro \clh_K$ by
\[
(X h) (w) = ( \langle h, \gamma_{1, w} \rangle_{\clh}, \ldots,
\langle h, \gamma_{m, w} \rangle_{\clh}) \in \mathbb{C}^m \quad
\quad (w \in \Omega).
\]
It follows that $X$ is a unitary and
\[
X M_z = M_z X.
\]
Also by the definition of $K$ we have
\[
K(\cdot, w) \eta = (\langle \gamma_{j, w}, \gamma_{i, z}
\rangle_{\clh} ) \eta,
\]
and
\[
\langle f(w), \eta \rangle_{\mathbb{C}^m} = \langle f, K(\cdot, w)
\eta \rangle_{\clh_K},
\]
for all $w \in \Omega, \eta \in \mathbb{C}^m$, and $f \in \clh_K$.
Moreover, the evaluation operator $ev_{w} : \clh_K \raro
\mathbb{C}^m$, $w \in \Omega$, satisfies
\[
\langle ev_w^* \eta, f \rangle_{\clh_K} = \langle \eta, f(w)
\rangle_{\mathbb{C}^m} = \langle K(\cdot, w) \eta, f
\rangle_{\clh_K},
\]
for $\eta \in \mathbb{C}^m$ and $f \in \clh_K$. In particular,
\[
ev_w^* \eta = K(\cdot, w) \eta,
\]
and hence
\[
K(z,w) = ev_z \circ ev_w^*,
\]
for all $z, w \in \Omega, \eta \in \mathbb{C}^m$. Therefore, we have
prove the following  (see \cite{Alpay}, \cite{CS}):

\begin{Theorem}
Let $\clh \in B_m(\Omega)$, and let $\{\gamma_{i, w}\}_{i=1}^m$ be
an anti-holomorphic frame of $\clh$. If
\[
K(z, w) = (\langle \gamma_{j, w}, \gamma_{i, z} \rangle_{\clh} )_{m
\times m} \quad \quad (z, w \in \Omega),
\]
then $\clh_K$ is a reproducing kernel Hilbert  module. Moreover,
$M_z$ on $\clh$ and $M_z$ on $\clh_K$ are unitarily equivalent.

\end{Theorem}

\textsf{By virtue of the above theorem, from now on, we will often
use the representation $\clh_K$ of $\clh$ in $B_m(\Omega)$.}

\newsection{Derivatives of Holomorphic Maps and Curvatures}

The purpose of this section is to study the curvatures, a projection
formula for eigenspace bundles and a trace-curvature formula in
terms of Hilbert Schmidth norm of derivatives of eigenspace bundles.
Most of the results of this section are known. However, our method
of proofs are more geometric and explicit.

Let $\clh \in B_m(\Omega)$, and let $\{\gamma_{i, w} : 1 \leq i \leq
m\}$ be an anti-holomorphic frame of $E_{\clh}$. The curvature
matrix of $E_{\clh}$ is given by
\[
\clk_{\clh}(z) = - \overline{\partial}[G^{-1}(\bar w) \partial G(\bar w)] \quad
\quad (z \in \Omega),
\]
where $G$ is the Gram matrix given by (see Section 2)
\[
G(w) = (\langle \gamma_{j,w}, \gamma_{i, w}\rangle_{\clh})_{i,j =
1}^m = K(w, w) \quad\quad (w \in \Omega).
\]
If $E_{\clh}$ is a line bundle then it follows that
\[
\clk_{\clh}(w) = - \frac{\partial^2}{\partial w \partial \bar{w}}
\log\|\gamma_w\|^2 \quad \quad (w \in \Omega).
\]

Let $\clh$ and $\clk$ be Hilbert spaces, and let $T : \clh \raro
\clk$ be a left invertible operator. Then
\[
L = (T^* T)^{-1} T^*,
\]
is a left inverse of $T$, and hence
\[
Q = T L,
\]
is an orthogonal projection of $\clk$ onto $\mbox{ran~} T$, that is,
\[
Q = P_{\mbox{ran~} T}.
\]

Now let $\clh \in B_m(\Omega)$. Then applying the above observations
to $\Gamma(w) = ev_w^* : \mathbb{C}^m \raro \clh_K$, $w \in \Omega$,
we have the following useful result.

\begin{Theorem} \textsf{(Projection Formula)}
Let $\clh \in B_m(\Omega)$, $\Gamma(w) = ev_w^*$, and let
\[
G(w) = \Gamma(w)^* \Gamma(w),
\]
for all $w \in \Omega$. Then
\[
P_{\mbox{ker} (M_z - w)^*} = \Gamma(w)  G(w)^{-1} \Gamma(w)^*,
\]
and
\[
\mbox{ker} (M_z - w)^* = \mbox{~ran~} ev_w^*.
\]
\end{Theorem}
\NI For simplicity of notation, we denote
\[
\Pi_{\clh}(w) = \Pi(w) = P_{\mbox{ker~}(M_z - w)^*} = \Gamma(w)
G(w)^{-1} \Gamma(w)^* \quad \quad (w \in \Omega).
\]

The above theorem should be compared with the results of Curto and
Salinas (see, for example, Theorem 2.2 in \cite{CS}).

The following result provides a useful relation between curvature
and the derivatives of the orthogonal projection-valued map $\Pi$.

\NI We first define a holomorphic map $\overline{\Pi}$ as follows
\[
\overline{\Pi}(w) = \Pi(\bar w) \quad \quad (w\in \Omega).
\]
Then
\[ \overline{\Pi}(w) = P_{\mbox{ker~}(M^*_z - w)} = \Gamma(\bar
w) G(\bar w)^{-1} \Gamma^*(\bar w)\quad \quad (w \in \Omega)
\]
Also define $\overline{\Gamma} $ as
\[
\overline{\Gamma}(w) = \Gamma(\bar w) \quad \quad (w\in \Omega).
\]
Then we have:
\begin{Theorem}\label{Prop DPi}
If $\clh \in B_m(\Omega)$, then
\[
\bar{\partial} \overline{\Pi} (w) \partial \overline{\Pi} (w) = -
\overline{\Gamma} (w) [\clk_{\clh}(w) G^{-1}(\bar w)]
\overline{\Gamma}^*(w) \quad \quad (w \in \Omega).
\]
\end{Theorem}
\NI\textsf{Proof.} Note that
\[
\partial \overline{\Pi} (w)= [(\partial \overline{\Gamma} ( w)) G^{-1}(\bar w) + \overline{\Gamma} ( w) (\partial
G^{-1})(\bar w)] \overline{\Gamma}^*( w),
\]
and so
\[
\begin{split}
(\bar{\partial} \overline{\Pi}) (w)(\partial \overline{\Pi})(w) & =
\overline{\Gamma}(w) [(\bar{\partial} G^{-1}(\bar w))
\overline{\Gamma}^* (w) + G^{-1}(\bar w) (\bar{\partial}
\overline{\Gamma}^*)(w)] [(\partial \overline{\Gamma}(w))
G^{-1}(\bar w) +
\\
& \;\;\; \; \overline{\Gamma} (w) (\partial G^{-1})(\bar w)]  \bar
\Gamma^*(w)
\\
& = \overline{\Gamma}(w) [(\bar{\p} G^{-1}(\bar w))
\overline{\Gamma}^*( w) (\p \overline{\Gamma}(w)) G^{-1} (\bar w)+
(\bar{\p} G^{-1}(\bar w)) \overline{\Gamma}^*(w) \overline{\Gamma}
(w)(\p G^{-1}(\bar w)) +
\\
& \;\;\; \; G^{-1}(\bar w) (\bar{\p} \overline{\Gamma}^*) (\p
\overline{\Gamma}(w)) G^{-1}(\bar w) + G^{-1} (\bar w)(\bar{\p}
\overline{\Gamma}^*(w)) \overline{\Gamma} (w) (\p G^{-1})(\bar w)]
\overline{\Gamma}^*(w)
\\
& = \overline{\Gamma}(w) [ (\bar{\p} G^{-1}(\bar w)) (\p G(\bar w))
G^{-1} (\bar w)+ (\bar{\p} G^{-1}(\bar w)) G (\bar w) (\p
G^{-1}(\bar w)) \\
&+ G^{-1}(\bar w) (\p \bar{\p} G(\bar w)) G^{-1} (\bar w)+
G^{-1}(\bar w) (\bar{\p} G(\bar w)) (\p G^{-1})(\bar w)] \G^*(\bar
w)
\\
& =\overline{\Gamma}(w) \{[ (\bar{\p} G^{-1}) (\bar w) (\p G) (\bar
w) + G^{-1} (\bar w) (\p \bar{\p} G) (\bar w)] G^{-1}(\bar w) +
\\
& \;\;\; \; [(\bar{\p} G^{-1}) (\bar w) G (\bar w) + G^{-1}(\bar w)
(\bar{\p} G)(\bar w)] (\p G^{-1})(\bar w) \} \overline{\Gamma}^*(w).
\end{split}
\]
Then
\[
(\bar{\partial} \overline{\Pi}) (w)(\partial \overline{\Pi})(w) =
\overline{\G}(w) [ (\bar{\p} G^{-1}) (\p G) + G^{-1} (\p \bar{\p}
G)] (\bar w)G^{-1} (\bar w)\overline{\G}^*(w),
\]
because
\[
(\bar{\p} G^{-1}) G + G^{-1} (\bar{\p} G) = \bar{\p}(G^{-1} G) =
\bar{\p}(I) = 0.
\]
Therefore
\[
\begin{split}
(\bar{\partial} \overline{\Pi})(w) (\partial \overline{\Pi})(w) & =
\overline{\Gamma}(w) [ (\bar{\p} G^{-1}) (\p G) + G^{-1} (\p
\bar{\p} G)](\bar w) G^{-1}(\bar w) \overline{\Gamma}^*(w)
\\
& = \overline{\Gamma} (w) [\bar{\p}(G^{-1} (\bar w)\p G(\bar w))]
G^{-1} (\bar w)\overline{\Gamma}^*(w)
\\
& = -\overline{\Gamma}(w) [\clk_{\clh} (w)G^{-1}(\bar w)]
\overline{\Gamma}^*(w).
\end{split}
\]
This completes the proof. \qed

As a consequence, we have the following equality:

\begin{Corollary}\label{HJK}
Let ${\mathcal H} \in B_m({\Omega})$. Then
\[
\|{ \partial}\overline{\Pi}(w)\|_{2}^2 = -\text{trace~
}\mathcal{K}_{\mathcal H}({w}),
\]
where $\|{\partial}\overline{\Pi}(w)\|_{2}^2$ is the Hilbert-Schmidt
norm of ${ \partial} \overline{\Pi}(w)$ and $w \in \Omega$.
\end{Corollary}
\NI\textsf{Proof.} Clearly $\clk_{\clh}(w)$, $w \in \Omega$, is a
finite rank operator on $\clh$. Hence, in particular,
$\clk_{\clh}(w)$, $w \in \Omega$, is in trace class. From this, it
follows that
\[
\begin{split}
\text{trace~} (\overline{\Gamma}(w) [\clk_{\clh} (w)G^{-1}(\bar w)]
\overline{\Gamma}^*(w)) & = \text{trace~} ([\clk_{\clh}
(w)G^{-1}(\bar w)] \overline{\Gamma}^*(w) \overline{\Gamma}(w))
\\
& = \text{trace~} \clk_{\clh}(w).
\end{split}
\]
Furthermore, Theorem \ref{Prop DPi} shows that
\[
\|{\partial}\overline{\Pi}(w)\|_{2}^2=\text{trace~}((\bar{\partial}
\overline{\Pi}) (\partial \overline{\Pi})) = - \text{trace~} (
(\overline{\Gamma}(w) [\clk_{\clh} (w)G^{-1}(\bar w)]
\overline{\Gamma}^*(w))).
\]
Thus
\[
\|{\partial}\overline{\Pi}(w)\|_{2}^2 = - \text{trace~}
\clk_{\clh}(w),
\]
for all $w \in \Omega$. This completes the proof. \qed

The above result is due to the first author and Hou and Kwon
\cite{HJK} (for the Hardy and the weighted Bergman spaces, see Lemma
1.7 in \cite{KT} and Lemma 3.3 in \cite{DKT}, respectively).
However, the present proof is more geometric and simple.

\newsection{Cowen-Douglas Atoms}

In this section, we introduce the concept of a Cowen-Douglas atom, a
special but large class of rank one Cowen-Douglas Hilbert modules
over $\mathbb{D}$. Moreover, a Cowen-Douglas atom admits
$\frac{1}{K}$-calculus in the sense of Arazy and Englis \cite{AE}.
Our presentation of $\frac{1}{K}$-calculus is restricted to one
variable. For more details, we refer the readers to the work by
Arazy and Englis \cite{AE}.

\begin{Definition}\label{atom-def}
A Hilbert module $\clm \in B_1(\mathbb{D})$ is said to be a
Cowen-Douglas atom if

(i) the set of polynomials $\mathbb{C}[z]$ is dense in $\clm$,

(ii) there exists a sequence of polynomials $\{p_l(z,\bar{w})\}_l
\subseteq \mathbb{C}[z, \bar{w}]$ such that
\[p_l(z,\bar{w}) \raro \frac{1}{k_{\clm}(z,w)},\]as $l \raro \infty$ and
for all $z, \bar{w} \in \mathbb{D}$, where $k_{\clm}$ is the kernel
function of $\clm$ (see Section 2),

(iii) $\mbox{sup}_l\|p_l(M_z, M_z^*)\| < \infty$, and

(iv) $\{M_z\}' = \{M_{\varphi} : \varphi \in
H^{\infty}(\mathbb{D})\}$.

\end{Definition}

Appealing to Theorem 1.6 of \cite{AE}, a Cowen-Douglas atom $\clm$
admits a $\frac{1}{K}$-calculus. Here we do not intend to define the
$\frac{1}{K}$-calculus but spell out the required properties of such
concept in the present set up. We again refer the reader to
\cite{AE} for details.

Note that by condition (i) in the above definition and the
Gram-Schmidth orthogonalization process, for a Cowen-Douglas atom
$\clm$ there exists a sequence of orthonormal basis of polynomials
$\{q_l(z) : l \geq 0\}$ such that
\begin{equation}\label{kernel-basis}
k_{\clm}(z, w) = \sum_{l \geq 0} q_l(z)
\overline{q_l(w)}.
\end{equation}

\textit{We henceforth assume $\clm$ to be a fixed Cowen-Douglas atom
with the sequence of polynomials $\{p_l(z,\bar{w})\}$ as in (ii) of
Definition \ref{atom-def} and the orthonormal basis as above with
the kernel function identity (\ref{kernel-basis}).}

Natural examples of Cowen-Douglas atoms include the Hardy space and
the weighted Bergman spaces (cf. \cite{AE}).

We turn now to define an analogue of the contractive Hilbert
modules.

\begin{Definition}\label{M-def}
A Hilbert module $\clh$ over $\mathbb{C}[z]$ is said to be
$\clm$-contractive if

(i) $\mbox{sup}_l\|p_l(M_z, M_z^*)\| < \infty$ and

(ii) $C_{\clh} := \mbox{WOT}-\lim_{l \raro \infty} p_l(M_z, M_z^*)$
is a positive operator.
\end{Definition}

The following lemma will be useful in the sequel.

\begin{Lemma}\label{C}
Let $\clm$ be a Cowen-Douglas atom, and let $\cle$ be a Hilbert
space. Also let $\clq = (\clm \otimes \cle)/\cls$ be an
$\clm$-contractive Hilbert module for some submodule $\cls$ of $\clm
\otimes \cle$. Then,

(i) $C_{\clm \otimes \cle} = I_{\clm} \otimes P_{\cle}$, and

(ii) $C_{\clq} = P_{\clq}C_{\clm \otimes \cle} P_{\clq} = P_{\clq}
(I_{\clm} \otimes P_{\cle})P_{\clq}$.
\end{Lemma}

\NI\textsf{Proof.} Let $z, w \in \mathbb{D}$, and let $x, y \in
\cle$. Then for all $l \geq 1$ we have
\[
\begin{split}
\langle p_l(M_z \otimes I_{\cle}, M_z^* \otimes I_{\cle})
(k_{\clm}(\cdot, w) \otimes x), k_{\clm}(\cdot, z) \otimes y\rangle
& = p_l(z, \bar{w}) \langle
k_{\clm}(\cdot, w) \otimes x, k_{\clm}(\cdot, z) \otimes y\rangle \\
& = p_l(z, \bar{w}) k_{\clm}(z, w) \langle x, y\rangle.
\end{split}
\]
This shows, by letting $l \raro \infty$, that
\[
\begin{split}
\langle C_{\clm \otimes \cle}(k_{\clm}(\cdot, w) \otimes x),
k_{\clm}(\cdot, z) \otimes y\rangle & = \frac{1}{k_{\clm}(z,
w)}k_{\clm}(z, w) \langle x, y\rangle = \langle x, y\rangle \\& =
\langle (I_{\clm} \otimes P_{\cle}) (k_{\clm}(\cdot, w) \otimes x),
k_{\clm}(\cdot, z) \otimes y\rangle,
\end{split}
\]
and hence $C_{\clm \otimes \cle} = I_{\clm} \otimes P_{\cle}$. This
completes the proof of part (i).

\NI To prove (ii) we compute \[p_l(P_{\clq}(M_z \otimes
I_{\cle})|_{\clq}, P_{\clq}(M_z^* \otimes I_{\cle})|_{\clq}) =
P_{\clq}(p_l(M_z \otimes I_{\cle}, M_z^* \otimes
I_{\cle}))|_{\clq}.\]Letting $l \raro \infty$ in WOT, we deduce from
part (i) that
\[
C_{\clq} = P_{\clq}(I \otimes P_{\cle})|_{\clq}.
\]
This completes the proof of the lemma.
\qed

We need the following analogue of $C_{\cdot 0}$ contractions
\cite{NF}.

Let $A$ be a positive operator on a Hilbert space $\clh$. Define
$Q_{l,A} : \clb(\clh) \raro \clb(\clh)$ for all $l \geq 1$ by
\[
Q_{l, A}(T) = I - \sum_{0 \leq j <l}q_j(T) A q_j(T)^* \quad\quad (T
\in \clb(\clh)).
\]
An $\clm$-contractive Hilbert module $\clh$ is said to be
\textit{pure} if
\[
\mbox{SOT}- \lim_{l \raro \infty} Q_{l,
C_{\clh}}(M_z) = 0.
\]

Let $\cle$ be a Hilbert space and let $\clh = \clm \otimes \cle$.
Observe that if $\clq$ is a quotient module of $\clh$ then $\clq$ is
a pure $\clm$-contractive Hilbert module (see \cite{AE}).

\newsection{Quotient modules and tensor product bundles}

The aim of the present section is to prove that the hermitian
anti-holomorphic vector bundle of a pure $\clm$-contractive Hilbert
module in $B_m(\mathbb{D})$ can be represented as the tensor product
bundle of a hermitian anti-holomorphic line bundle and a rank $m$
hermitian anti-holomorphic vector bundle.

We start by recalling a version of the model theorem due to Arazy
and Englis (Corollary 3.2 in \cite{AE}).

\begin{Theorem}\label{A-E}(\textsf{Arazy-Englis})
Let $\clh$ be a Hilbert module over $\mathbb{C}[z]$. Then $\clh$ is
a pure $\clm$-contractive Hilbert module if and only if $\clh$ is
unitarily equivalent to a quotient module of $\clm \otimes \cle$ for
some Hilbert space $\cle$.
\end{Theorem}

The following proposition shows that an $\clm$-contractive Hilbert
modules in $B_m(\mathbb{D})$ is pure.

\begin{Proposition}\label{pure}
Let $\clh \in B_m(\mathbb{D})$ be an $\clm$-contractive Hilbert
module. Then $\clh$ is pure.
\end{Proposition}

\NI\textsf{Proof.} Let $\{\gamma_{i, w} : 1 \leq i \leq m\}$ be an
anti-holomorphic frame of $E_{\clh}$ with
\[
M_z^* \gamma_{i, w} = \bar{w} \gamma_{i, w},
\]
for all $w \in \mathbb{D}$ and $1 \leq i \leq m$. Then for all $z, w
\in \mathbb{D}$, and $1 \leq i, j \leq m$, we have
\[
\begin{split}
\langle C_{\clh} \gamma_{i, w}, \gamma_{j, z} \rangle & = \lim_{l
\raro \infty} \langle p_l(M_z, M_z^*) \gamma_{i, w}, \gamma_{j, z}
\rangle
\\
& = (\lim_{l \raro \infty} p_l(z, \bar{w})) \langle \gamma_{i, w},
\gamma_{j, z} \rangle
\\
& = \frac{1}{k_{\clm}(z, w)} \langle \gamma_{i, w}, \gamma_{j, z}
\rangle,
\end{split}
\]
and hence
\[
\begin{split}
\langle Q_{l, C_{\clh}}(M_z) \gamma_{i, w},
\gamma_{j, z}\rangle & = \langle \gamma_{i, w}, \gamma_{j, z}
\rangle - \langle \sum_{t=0}^{l-1} q_t(M_z) C_{\clh} q_t(M_z)^*
\gamma_{i, w}, \gamma_{j, z}\rangle \\& = \langle \gamma_{i, w},
\gamma_{j, z} \rangle - \langle \sum_{t=0}^{l-1} C_{\clh}
\overline{q_t(w)} \gamma_{i, w}, \overline{q_t(z)} \gamma_{j,
z}\rangle
\\& = \langle \gamma_{i, w}, \gamma_{j, z} \rangle -
(\sum_{t=0}^{l-1} q_t(z) \overline{q_t(w)}) \langle C_{\clh}
\gamma_{i, w}, \gamma_{j, z}\rangle \\& = \langle \gamma_{i, w},
\gamma_{j, z} \rangle - (\sum_{t=0}^{l-1} q_t(z) \overline{q_t(w)})
\frac{1}{k_{\clm}(z, w)} \langle \gamma_{i, w}, \gamma_{j, z}\rangle
\\& = (1 - (\sum_{t=0}^{l-1} q_t(z) \overline{q_t(w)}) \frac{1}{k_{\clm}(z,
w)})\langle \gamma_{i, w}, \gamma_{j, z}\rangle
\\& \raro 0 \mbox{~as~} l \raro \infty.
\end{split}
\]
From this we deduce that $Q_{l, C_{\clh}}(M_z) \raro 0$ in SOT.
This concludes the proof. \qed

Let $\clh \in B_m(\mathbb{D})$ be an $\clm$-contractive module. As
an application of the previous proposition and Theorem \ref{A-E},
$\clh$ can be realized as
\[
\clh \cong \clq := (\clm \otimes \cle)/ \cls,
\]
for some Hilbert space $\cle$ and submodule $\cls$ of $\clm \otimes
\cle$. That is,
\[
0 \raro \cls \raro \clm \otimes \cle \raro \clh \raro 0.
\]
Therefore, an $\clm$-contractive Hilbert module $\clh \in
B_m(\mathbb{D})$ can be realized as a quotient module $\clq$ of
$\clm \otimes \cle$ for some coefficient space $\cle$. In this
representation, the module map $M_z$ on $\clh$ is identified with
the compressed multiplication operator $P_{\clq}(M_z \otimes
I_{\cle})|_{\clq}$. Moreover,
\[P_{\clq} (M_z \otimes I_{\cle})^*|_{\clq} = (M_z \otimes
I_{\cle})^*|_{\clq}.\]\textit{In the rest of this paper we will
assume the quotient module representations of the class of pure
$\clm$-contractive Hilbert modules in $B_m(\mathbb{D})$.}

Also, often we will identify a Cowen-Douglas atom $\clm$ with the
reproducing kernel Hilbert module $\clh_{k_{\clm}}$ (see Section 2)
with section
\[
w \mapsto k_{\clm}(\cdot, w) \quad \quad (w \in \mathbb{D}).
\]

Now we are in a position to prove the main result of this section.

\begin{Theorem}\label{tensor-bundle}
Let $\cle$ be a Hilbert space, and let $\clq$ is a quotient module
of $\clm \otimes \cle$. Then $\clq \in B_m(\mathbb{D})$ if and only
if there exists a rank $m$ hermitian anti-holomorphic vector bundle
$V$ over $\mathbb{D}$ such that
\[
E_{\clq} \cong E_{\clm} \otimes V.
\]
Moreover, in this case
\[
\clk_{\clq} = \clk_{\clm} + \clk_V .
\]
\end{Theorem}

\NI\textsf{Proof.} The sufficient part of the statement is trivial,
so we only have to prove the necessary part. Let $\{\gamma_{i, w} :
1 \leq i \leq m\}$ be an anti-holomorphic frame of $E_{\clq}$ such
that \[M_z^* \gamma_{i, w} = \bar{w} \gamma_{i, w},\]for all $1 \leq
i \leq m$ and $w \in \mathbb{D}$. Then for all $l \geq 1$ we have
\[p_l(M_z, M_z^*) \gamma_{i, w} = p_l(z, \bar{w}) \gamma_{i, w}.\]
Letting $l \raro \infty$ in WOT, and applying Lemma \ref{C} we have
\[\frac{1}{k_{\clm}(\cdot,w)} \gamma_{i,w} = C_{\clq} \gamma_{i,w} =
P_{\clq}(I \otimes P_{\cle}) \gamma_{i,w}.\]Since \[ P_{\clq}(I
\otimes P_{\cle}) \gamma_{i,w} = \gamma_{i,w}(0),\]we have
\[\frac{1}{k_{\clm}(\cdot,w)} \gamma_{i,w} = \gamma_{i,w}(0).\]
Therefore
\begin{equation}\label{C-section} \gamma_{i, w} = k_{\clm}(\cdot, w)
\otimes \gamma_{i, w}(0) = k_{\clm}(\cdot, w) \otimes v_{i, w},
\end{equation}where $v_{i, w} := \gamma_{i, w}(0)$ for all $1 \leq i
\leq m$ and $w \in \mathbb{D}$. Let $V$ be the anti-holomorphic
curve over $\mathbb{D}$ with
\[
V(w) = \mbox{span~} \{v_{i,w} : 1 \leq i \leq m\} \quad \quad (w \in
\mathbb{D}).
\]
Then we conclude that $E_{\clq} \cong E_{\clm} \otimes V$.

\NI Finally, let $G_V$ be the Gram matrix corresponding to the frame
$\{v_{i,w}\}$ of $E_V$. Then
\[
\begin{split}
\clk_{E_{\clq}}(w) & = - \overline{\partial}[G_{E_{\clm}}^{-1}(\bar w)
\partial G_{E_{\clm}}(\bar w)]
\\
& = - \overline{\partial}[\frac{1}{\|k_{\clm}(\cdot, w)\|^2}
G_{V}^{-1}(\bar w)
\partial \{\|k_{\clm}(\cdot, w)\|^2 G_{V}(\bar w)\}]\\& = -
\overline{\partial}[\frac{1}{\|k_{\clm}(\cdot, w)\|^2}
G_{V}^{-1}(\bar w)\{\partial(\|k_{\clm}(\cdot, w)\|^2) G_V(\bar w) +
\|k_{\clm}(\cdot, w)\|^2 \partial G_V(\bar w)\}] \\ & = -
\overline{\partial}[\frac{1}{\|k_{\clm}(\cdot, w)\|^2}
\partial(\|k_{\clm}(\cdot, w)\|^2) + G_{V}^{-1}(\bar w) \partial(G_V(\bar w))] \\ &
= \clk_{E_{\clm}} (w) + \clk_{V}(w),
\end{split}
\]
for all $w \in \mathbb{D}$. This concludes the proof of the theorem.
\qed

In particular, we have the following useful result.

\begin{Corollary}\label{C2}
Let $\clh \in B_m(\mathbb{D})$ be a pure $\clm$-contractive Hilbert
module. Then there exists a hermitian anti-holomorphic vector bundle
$V$ of rank $m$ over $\mathbb{D}$ such that
\[
E_{\clh} \cong E_{\clm} \otimes V.
\]
Moreover
\[
\Pi_{\clh} = \Pi_{\clm}\otimes \Pi_{V},
\]
and for each $w \in \mathbb{D}$,
\[
\|{\partial}\Pi_{\clh}(w)\|^2_2 = m \|{\partial}\Pi_{\clm}(w)\|^2_2
+ \|{\partial}\Pi_V(w)\|^2_2 = m |\clk_{\clm}(w)| +
\|{\partial}\Pi_V(w)\|^2_2.
\]
\end{Corollary}

\NI\textsf{Proof.} First two conclusions follows from Theorem
\ref{tensor-bundle}. For the remaining parts, we follow the same
line of arguments as in \cite{KT} or \cite{DKT}. Since
\[
{\partial}\Pi_{\clq}(w) =
{\partial}(\Pi_{\clm}(w) \otimes \Pi_V(w)) = {\partial}\Pi_{\clm}(w)
\otimes \Pi_V(w) + \Pi_{\clm}(w) \otimes {\partial}\Pi_{V}(w),
\]
we have that
\[
\begin{split}
\|{\partial}\Pi_{\clq}(w)\|^2_2  = & tr([{\partial}\Pi_{\clm}(w)
\otimes \Pi_V(w)][{\partial}\Pi_{\clm}(w) \otimes \Pi_V(w)]^*) + 2
\mbox{Real}\{ tr([{\partial}\Pi_{\clm}(w) \otimes \Pi_V(w)]^*\\&
[\Pi_{\clm}(w) \otimes {\partial}\Pi_{V}(w)])\} + tr([\Pi_{\clm}(w)
\otimes {\partial}\Pi_{V}(w)]^* [\Pi_{\clm}(w) \otimes
{\partial}\Pi_{V}(w)]).
\end{split}
\]
Notice that $\overline{\partial}\Pi_{\clm}(w)
\Pi_{\clm}(w) = 0$ and hence the middle term in the last expression
vanishes. Therefore,
\[
\begin{split}
\|{\partial}\Pi_{\clq}(w)\|^2_2 & = \|{\partial}\Pi_{\clm}(w)
\otimes \Pi_V(w)\|^2_2 + \|\Pi_{\clm}(w) \otimes
{\partial}\Pi_{V}(w)\|^2_2\\& = \|{\partial}\Pi_{\clm}(w)\|^2_2
\|\Pi_V(w)\|^2_2 + \|\Pi_{\clm}(w)\|^2_2
\|{\partial}\Pi_{V}(w)\|^2_2 \\& = m |\clk_{\clm}(w)| +
\|{\partial}\Pi_V(w)\|^2_2,
\end{split}
\]
where the last equality follows from Corollary \ref{HJK}. This
completes the proof. \qed

\newsection{Quasi-affinity}

In this section we discuss the issue of quasi-affinity of Hilbert
modules in the Cowen-Douglas class $B_m(\mathbb{D})$. We begin with
the definition of quasi-affinity.

Let $\clh$ and $\clk$ be two Hilbert modules. Then we say that
$\clh$ is \textit{quasi-affine} to $\clk$, and denote by $\clh \prec
\clk$, if there exists a module map $X : \clh \raro \clk$ such that
$X$ is one-to-one and has dense range.

\begin{Theorem}\label{q-s}
Let $\clh \in B_m(\mathbb{D})$ be a pure $\clm$-contractive Hilbert
module and $\{\gamma_{i,w}\}_{i=1}^m$ be an anti-holomorphic frame
of $E_{\clh}$ such that \[\sup_{w \in \mathbb{D}}
(\frac{\|\gamma_{i,w}\|}{\|k_{\clm}(\cdot, w)\|}) < \infty,\]for all
$i = 1, \ldots, m$. Then

(i) there exists a one-to-one module map $X : \clh \raro \clm
\otimes \mathbb{C}^m$, and

(ii) $\clh \prec \cls$ for some submodule $\cls \subseteq \clm
\otimes \mathbb{C}^m$.
\end{Theorem}

\NI\textsf{Proof.} Identifying $\clh$ with $\clq = (\clm \otimes
\cle)/\cls$ for some submodule $\cls$ of $\clm \otimes \cle$, we let
\[
\gamma_{i,w} = k_{\clm}(\cdot, w) \otimes v_{i,w},
\]
for each $i = 1, \ldots, m$. Set
\[
\delta: = \sup_{w \in \mathbb{D}}
(\frac{\|\gamma_{i,w}\|}{\|k_{\clm}(\cdot, w)\|}) = \sup_{w \in
\mathbb{D}} \|v_{i,w}\| < \infty.
\]
For each $z \in \mathbb{D}$,
define $\Theta(z) \in \clb(\cle, \mathbb{C}^m)$ by
\[
\Theta(z) \eta = (\langle \eta, v_{1,z}\rangle_{\cle}, \ldots,
\langle \eta, v_{m,z}\rangle_{\cle}) \in \mathbb{C}^m.
\]
Then
\[
\|\Theta(z) \eta\|^2 = \sum_{i=1}^m |\langle \eta,
v_{i,z}\rangle_{\cle}|^2 \leq \|\eta\|^2 \sum_{i=1}^m \|v_{i,w}\|^2
\leq m \delta^2 \|\eta\|^2,
\]
for all $\eta \in \cle$. Consequently, $\Theta \in
H^{\infty}_{\clb(\cle, \mathbb{C}^m)}(\mathbb{D})$. Furthermore, for
$f \in \cls = \clq^{\perp}$ and $w \in \mathbb{D}$ we have
\[
\begin{split}(\Theta f)(w) & = \Theta(w) f(w)= (\langle f(w),
v_{1,w}\rangle_{\cle}, \ldots, \langle f(w), v_{m,w}\rangle_{\cle})
\\& = (\langle f, k_{\clm}(\cdot, w) \otimes v_{1,w}\rangle_{\clm \otimes \cle},
\ldots, \langle f, k_{\clm}(\cdot, w) \otimes v_{m,w}\rangle_{\clm
\otimes \cle}) \\& = (\langle f, \gamma_{1,w}\rangle_{\clm \otimes
\cle}, \ldots, \langle f, \gamma_{m,w}\rangle_{\clm \otimes
\cle})\\& = 0.\end{split}
\]
Hence, $M_{\Theta}\cls = \{0\}$. Next we
define $X : \clq \raro \clm \otimes \mathbb{C}^m$ by \[Xf =
M_{\Theta}f,\] for all $f \in \clq$. Then
\[
\begin{split}
X P_{\clq}(M_z \otimes I_{\cle})|_{\clq} & = M_{\Theta} P_{\clq}(M_z
\otimes I_{\cle})|_{\clq}
\\
& = M_{\Theta} (M_z \otimes I_{\cle})|_{\clq}
\\
& = (M_z \otimes I_{\mathbb{C}^m}) M_{\Theta}|_{\clq},
\end{split}
\]
that is,
\[
\begin{split}
X P_{\clq}(M_z \otimes I_{\cle})|_{\clq} = (M_z \otimes
I_{\mathbb{C}^m})X,
\end{split}
\]
and hence, $X$ is a module map. To prove that $X$ is one-to-one, or
equivalently, that $X^*$ has dense range, we compute
\[
\langle \Theta(w)^* e_i, \eta\rangle_{\cle} = \langle e_i, \Theta(w)
\eta\rangle_{\mathbb{C}^m} = \langle e_i, \sum_{j=1}^m \langle \eta,
v_{j,w}\rangle_{\cle} \rangle_{\mathbb{C}^m} = \langle v_{i,w},
\eta\rangle_{\cle},
\]
for all $w \in \mathbb{D}$, $\eta \in \cle$ and $i = 1, \ldots, m$.
Therefore
\[
\Theta(w)^* e_i = v_{i,w},
\]
and hence
\[
\begin{split}
X^*(k_{\clm}(\cdot, w) \otimes e_i) & = P_{\clq} M^*_{\Theta} (k_{\clm}(\cdot, w) \otimes
e_i)
\\
& = P_{\clq}(k_{\clm}(\cdot, w) \otimes \Theta(w)^* e_i)
\\
& = (k_{\clm}(\cdot, w) \otimes v_{i,w})
\\
& = P_{\clq} \gamma_{i,w}
\\
& = \gamma_{i,w}.
\end{split}
\]
Hence, $X$ is one-to-one. This proves part (i).

\NI Part (ii) follows from part (i) and by considering $\cls$ as the
range closure of $X$. \qed

In the anti-holomorphic vector bundle language, the above result can
be stated as follows : Suppose there exist an anti-holomorphic
bundle map $\Phi : E_{\clm \otimes \mathbb{C}^m} \raro E_{\clh}$ and
$\delta > 0$ such that
\[
\|\Phi(w) \eta_w\|_{\clh} \leq \delta \|\eta_w\|_{E_{\clm
\otimes \mathbb{C}^m}(w)},
\]
for all $\eta_w \in E_{\clm \otimes \mathbb{C}^m}(w)$ and $w \in
\mathbb{D}$. Then $\clh$ is quasi-affine to a submodule of $\clm
\otimes \cle$.

One might expect that the submodule $\cls$ in the above result is
the entire free module $\clm \otimes \mathbb{C}^m$. However, such
results are closely related with the issue of the
Beurling-Lax-Halmos type theorem for the Cowen-Douglas atoms. In
particular, for $\clm = H^2(\mathbb{D})$ the submodule $\cls$ is
unitarily equivalent with the Hardy module with the same
multiplicity as the rank of the map $\Theta(w)$ which is $m$.
Consequently, the conclusion is stronger for any
$H^2(\mathbb{D})$-contractive module, that is, $\clh$ is
quasi-affine to the Hardy module $H^2(\mathbb{D}) \otimes
\mathbb{C}^m$ (see \cite{U}). We point out that even the Bergman
module is quite subtle \cite{BCP} for this consideration.

\newsection{Similarity}

The purpose of this section is to relate the similarity problem with
the curvatures of Cowen-Douglas Hilbert modules. Recall that a
Hilbert module $\clh_1$ is said to be similar to a Hilbert module
$\clh_2$, denoted by $\clh_1\sim_s \clh_2$, if there exists an
invertible module map $X$ from $\clh_1$ to $\clh_2$.

We begin by generalizing a result by Uchiyama on similarity of a
contractive Hilbert modules to the Hardy module of finite
multiplicity (see Theorem 3.8 in \cite{U}).

\begin{Theorem}\label{sim-bundle}
Let $\clh \in B_m(\mathbb{D})$ be an $\clm$-contractive Hilbert
module. Then $\clh$ is similar to $\clm \otimes \mathbb{C}^m$ if and
only if there exists an anti-holomorphic pointwise invertible bundle
map $\Phi : E_{\clm \otimes \mathbb{C}^m} \raro E_{\clh}$ and
$\delta > 0$ such that
\[\frac{1}{\delta}\|\eta_w\|_{E_{\clm \otimes \mathbb{C}^m}(w)}
\leq \|\Phi(w) \eta_w\|_{\clh} \leq \delta \|\eta_w\|_{E_{\clm
\otimes \mathbb{C}^m}(w)},\] for all $\eta_w \in E_{\clm \otimes
\mathbb{C}^m}(w)$ and $w \in \mathbb{D}$.
\end{Theorem}

\NI\textsf{Proof.} Let $X : \clh \raro \clm \otimes \mathbb{C}^m$ be
an invertible module map. Then $\gamma_{i,w} := X^*(k_{\clm}(\cdot,
w) \otimes e_i)$ is the required anti-holomorphic frame of
$E_{\clh}$.

\NI For the converse, we proceed as in the proof of Theorem
\ref{q-s}. We first, consider an anti-holomorphic frame
$\{\gamma_{i,w}\}_{i=1}^m = \{k_{\clm}(\cdot, w) \otimes
v_{i,w}\}_{i=1}^m$ of $E_{\clh}$ and define $\Theta \in
H^{\infty}_{\clb(\cle, \mathbb{C}^m)}(\mathbb{D})$ by
\[\Theta(w) \eta = (\langle \eta, v_{1,w}\rangle_{\cle}, \ldots,
\langle \eta, v_{m,w}\rangle_{\cle}),\]for all $\eta \in \cle$ and
$w \in \mathbb{D}$. Now
\[
\|\Theta(w)^* x\| = \|\sum_{i=1}^m x_i
v_{i,w}\| = \frac{1}{\|k_{\clm}(\cdot, w)\|} \|\sum_{i=1}^m x_i
\gamma_{i,w}\|,
\]
and hence
\[
\|\Theta(w)^* x\| \geq \delta
\|x\|,
\]
for all $x \in \mathbb{C}^m$ and $w \in \mathbb{D}$. Hence $\Theta$
is right invertible (cf. Proposition 3.7 in \cite{U}). In
particular,
\[\mbox{ran} M_{\Theta} = \clm \otimes \mathbb{C}^m,\] and since
\[M_{\Theta} \cls = \{0\},\] the module map $X : \clq \raro \clm
\otimes \mathbb{C}^m$ defined by
\[
Xf = \Theta f \quad \quad (f \in \clq),
\]
is the required similarity. \qed

Now, we are ready to formulate the following similarity result for
pure $\clm$-contractive Hilbert modules. Applying our result to the
case where $\clh$ is the Hardy module, or a weighted Bergman module,
we recover the results of Kwon and Treil \cite{KT}, and Kwon, Treil
and Douglas \cite{DKT}.

\begin{Theorem}\label{similar-link}
Let $\clh \in B_m(\mathbb{D})$ be a pure $\clm$-contractive Hilbert
module. Then the following statements are equivalent:

(i) $\clh \sim_s \clm \otimes \mathbb{C}^m$.

(ii) There exists an anti-holomorphic pointwise invertible bundle
map $\Phi : E_{\clm \otimes \mathbb{C}^m} \raro E_{\clh}$ and
$\delta > 0$ such that \[\frac{1}{\delta}\|\eta_w\|_{E_{\clm \otimes
\mathbb{C}^m}(w)}  \leq \|\Phi(w) \eta_w\|_{\clh} \leq \delta
\|\eta_w\|_{E_{\clm \otimes \mathbb{C}^m}(w)},\] for all $\eta_w \in
E_{\clm \otimes \mathbb{C}^m}(w)$ and $w \in \mathbb{D}$.

(iii) There exists a bounded solution $\vp$ defined on $\mathbb{D}$
to the Poisson equation
\[
\Delta \vp = trace{\mathcal K}_{{\mathcal M}\otimes {\mathbb C}^m}-trace{\mathcal K}_{\mathcal H}.
\]
\end{Theorem}

\NI\textsf{Proof.}  The equivalence of (i) and (ii) is Theorem
\ref{sim-bundle}.

\NI (ii) implies (iii): We note that
\[
E_{\clm \otimes \mathbb{C}^m}(w) = \mbox{ker} (M_z - w)^* =
k_{\clm}(\cdot, w) \otimes \mathbb{C}^m,
\]
and
\[
E_{\clh}(w) = \mbox{ker}(M_z - w)^* = k_{\clm}(\cdot, w) \otimes
V(w).
\]
Consequently, for a given bundle equivalence $\Phi$ from $E_{\clm
\otimes \mathbb{C}^m}$ to $E_{\clh}$ there exists a one-to-one
bounded anti-holomorphic map $\Gamma : \mathbb{D} \raro
\clb(\mathbb{C}^m, \cle)$ such that
\[
\Phi(w) (k_{\clm}(\cdot, w) \otimes \eta) =
k_{\clm}(\cdot, w) \otimes \Gamma(w) \eta,
\]
or, equivalently,
\[
\Phi^{-1}(w) (k_{\clm}(\cdot, w) \otimes \Gamma(w) \eta)=
k_{\clm}(\cdot, w) \otimes \eta ,
\]
and
\[
V(w) = \mbox{ran}\Gamma(w),
\]
for all $\eta \in \mathbb{C}^m$ and $w \in \mathbb{D}$. Set
\[
F(w)=\Gamma(\overline{w}).
\]
Since
\[
\frac{1}{\delta}\|\eta_w\|_{E_{\clm \otimes \mathbb{C}^m}(w)}  \leq
\|\Phi(w) \eta_w\|_{\clh} \leq \delta \|\eta_w\|_{E_{\clm \otimes
\mathbb{C}^m}(w)} \quad \quad (w\in \mathbb{D}),
\]
by letting
\[
\eta_w=k_{\clm}(\cdot, w) \otimes \eta  \quad \quad (w\in
\mathbb{D}),
\]
we have
\[
\frac{1}{\delta}\|k_{\clm}(\cdot, w) \otimes
\eta\|_{E_{\clm \otimes
\mathbb{C}^m}(w)}  \leq \|k_{\clm}(\cdot, w) \otimes
\Gamma(w) \eta\|_{\clh} \leq \delta
\|k_{\clm}(\cdot, w) \otimes
\eta\|_{E_{\clm \otimes \mathbb{C}^m}(w)},
\]
that is
\[
\frac{1}{\delta^2}\|k_{\clm}(\cdot, w)\|^2 \|\eta\|^2 \leq
\|k_{\clm}(\cdot, w) \|^2 \|\Gamma(w) \eta\|^2_{\cle} \leq
\delta^2\|k_{\clm}(\cdot, w)\|^2 \|\eta\|^2,
\]
and so
\[
\frac{1}{\delta^2} \|\eta\|^2 \leq \langle \Gamma^*(w)\Gamma(w)
\eta, \eta\rangle \leq \delta^2 \|\eta\|^2,
\]
for all $\eta\in \mathbb{C}^m$. Writing $c=\delta^2$, the above
inequalities yields
\begin{equation}\label{ineq}
c^{-1}I\leq F^*F\leq cI.
\end{equation}

\textsf{Claim:} \textit{Let $\overline{\Pi}_{V}(w)$ be the
orthogonal projection of $\cle$ onto $V(\overline{w})$. Then $\|\p
\overline{\Pi}_{V}(w)\|_2\leq c^{\frac{1}{2}}\|F^{\prime}(w)\|_2$.}

\NI Indeed, as
\[
\overline{\Pi}_{V}=F(F^*F)^{-1}F^*,
\]
we have
\[
\overline{\Pi}_{V}F=F(F^*F)^{-1}F^*F=F.
\]
Then by a direct calculation, we have that
\begin{equation}\label{ineq1}
\partial\overline{\Pi}_{V}F=(I-\overline{\Pi}_{V})F^{\prime},
\end{equation}
and
\[
\partial \overline{\Pi}_{V}\overline{\Pi}_{V}=
\partial \overline{\Pi}_{V}.
\]
This yields
\[
\begin{array}{lll}
(I-\overline{\Pi}_{V})F^{\prime}(F^*F)^{-1}F^*&=&\partial \overline{\Pi}_{V}F(F^*F)^{-1}F^*\\
&=&\partial \overline{\Pi}_{V}\overline{\Pi}_{V}\\
&=&\partial \overline{\Pi}_{V}.
\end{array}
\]
By (\ref{ineq}) we have
\[
\begin{array}{llll}
\|\partial \overline{\Pi}_{V}\|_2&=&\|(I-\overline{\Pi}_{V}) F^{\prime}(F^*F)^{-1}F^*\|_2\\
&\leq &\|I-\overline{\Pi}_{V}\|\cdot \|F^{\prime}(F^*F)^{-1}F^*\|_2\\
&\leq &\|F^{\prime}(F^*F)^{-1}F^*\|_2\\
&\leq&\|(F^*F)^{-1}F^*\|\cdot \|F^{\prime}\|_2\\
&=&\|(F^*F)^{-1}F^*F(F^*F)^{-1}\|^{\frac{1}{2}}\cdot \|F^{\prime}\|_2\\
&=&\|(F^*F)^{-1}\|^{\frac{1}{2}}\cdot\|F^{\prime}\|_2\\
&\leq&c^{\frac{1}{2}}\|F^{\prime}\|_2.
\end{array}
\]
Therefore the claim does hold, as required.

\NI Finally, by Corollaries \ref{HJK} and \ref{C2}, we have
\[
trace{\mathcal K}_{{\mathcal M}\otimes {\mathbb
C}^m}(w)-trace{\mathcal K}_{\mathcal H}(w)=\|\partial
\overline{\Pi}_{V}(w)\|_2^2\leq c\|F^{\prime}(w)\|^2_2.
\]
Set
\[
\vp_1(w) = \|cF(w)\|_2^2.
\]
It follows that
\[
\Delta \vp_1(w) =c\|F^{\prime}(w)\|^2_2,
\]
and hence
\[
\Delta \vp_1 \geq trace{\mathcal K}_{{\mathcal M}\otimes {\mathbb
C}^m}-trace{\mathcal K}_{\mathcal H}.
\]
Let
\[
\clg_f(\lambda):= \frac{2}{\pi} \iint_{\mathbb{D}} \ln|\frac{w -
\lambda}{1 - \bar{\lambda} w}|f(w)dx dy,
\]
be the Green potential to the solution of $\Delta u=f(\lambda)$.
Then
\[
\vp_1=\clg_{\Delta \vp}\leq \clg_{trace{\mathcal K}_{{\mathcal
M}\otimes {\mathbb C}^m}-trace{\mathcal K}_{\mathcal H}}\leq 0.
\]
Set
\[
\vp=\clg_{trace{\mathcal K}_{{\mathcal M}\otimes {\mathbb
C}^m}-trace{\mathcal K}_{\mathcal H}}.
\]
Then $\vp$ is bounded and
\[
\Delta \vp=trace{\mathcal K}_{{\mathcal M}\otimes {\mathbb
C}^m}-trace{\mathcal K}_{\mathcal H}.
\]

\NI (iii) implies (i): We use Theorem 0.2 in \cite{TW} to get a
bounded anti-holomorphic projection $\Theta(w)$ onto $\mbox{ran}
\Pi_V(w)$. Let $\Theta_i$ be the inner part of the inner-outer
factorization of $\Theta$. Then it follows that the Toeplitz
operator $T_{\Theta_i}$ is invertible and the required similarity
operator with (see \cite{DKT} or \cite{KT} for more details)
\[
T_{\Theta_i}E_{\clm \otimes \mathbb{C}^m}(w)=E_{\clh}(w) ,  w\in
\mathbb{D}.
\]
\qed

It is of interest to note the following consequence of Theorem
\ref{similar-link}:

\NI Let $\clm$ and $\tilde{\clm}$ be two Cowen-Douglas atoms, and
let $\clh, \tilde{\clh} \in B_m(\mathbb{D})$ be $\clm$-contractive
and $\tilde{\clm}$ Hilbert modules, respectively. Let $V$ and
$\tilde{V}$ be the corresponding hermitian anti-holomorphic vector
bundles such that $E_{\clh} \cong E_{\clm} \otimes V$ and
$E_{\tilde{\clh}} \cong E_{\tilde{\clm}} \otimes \tilde{V}$ (see
Theorem \ref{tensor-bundle}).

\NI Now if $\clh$ is similar to $\clm \otimes \mathbb{C}^m$ then by
Corollary \ref{C2} and part (iii) of the previous theorem, we have
\[
\Delta \varphi (w) =  \| \overline{\partial} \Pi_V(w)\|_2^2,
\]
for some bounded subharmonic function on $\mathbb{D}$. Another
application of Corollary \ref{C2} and part (iii) of the previous
theorem to $\tilde{\clh}$ yields the similarity of $\tilde{\clh}$ to
$\tilde{\clm} \otimes \mathbb{C}^m$. Therefore, we have the
following result.

\begin{Corollary}
Let $\clm$ and $\tilde{\clm}$ be two Cowen-Douglas atoms, and let
$V$ be a rank $m$ hermitian anti-holomorphic vector bundle over
$\mathbb{D}$. Then the pure $\clm$-contractive Hilbert module $\clh$
corresponding to the hermitian anti-holomorphic vector bundle
$E_{\clm} \otimes V$ is similar to $\clm \otimes \mathbb{C}^m$ if
and only if the $\tilde{\clm}$-contractive Hilbert module
$\tilde{\clh}$ is similar to $\tilde{\clm} \otimes \mathbb{C}^m$,
where $\tilde{\clh} \in B_m(\mathbb{D})$ is the Hilbert module
corresponding to the hermitian anti-holomorphic vector bundle
$E_{\tilde{\clm}} \otimes V$.
\end{Corollary}

The above result is a generalization of Corollary 4.5 (restricted to
the Cowen-Douglas atoms) in \cite{DKKS2} where the quotient module
representations are assumed to be the orthocomplements of the
submodules implemented by left invertible multipliers. Moreover, the
free modules associated to the quotient modules are also assumed to
be of finite rank.

Let $\clh \in B_m(\mathbb{D})$ be a pure $\clm$-contractive Hilbert
module such that $E_{\clh} \cong E_{\clm} \otimes V$ and $V(w)
\subseteq \cle$, $w \in \mathbb{D}$. Let $\overline{\Pi}_V$(w) be
the orthogonal projection of $\cle$ onto $V(\bar{w})$. In the
following theorem, we prove that if $\overline{\Pi}_V$(w) factors as
\[
\overline{\Pi}_{V}(w) = F(w)G(w) \quad \quad (w\in {\mathbb{D}}),
\]
for some $F \in H^{\infty}_{B(\cle_*, \cle)}({\mathbb{D}})$ and  $G
\in L^{\infty}_{B(\cle, \cle_*)}({\mathbb{D}})$ with $ ran~F(w) =
V(w)$, $w \in \mathbb{D}$, then $\clh$ is similar to $\clm\otimes
{\mathbb C}^m$.

\begin{Theorem}\label{last}
Let $\cle$ and $\cle_*$ be two Hilbert spaces, let $\clh \in
B_m(\mathbb{D})$ be a pure $\clm$-contractive Hilbert module, and
let $F \in H^{\infty}_{B(\cle_*, \cle)}({\mathbb{D}})$. Let
\[
E_{\clh} \cong E_{\clm} \otimes V,
\]
where $V(w) \subseteq \cle$, and
\[
\mbox{ran~} F(w) = V(w) \quad \quad (w \in \mathbb{D}).
\]
Assume moreover that $\overline{\Pi}_V$(w) is the orthogonal
projection of $\cle$ onto $V(\bar{w})$, and let
\[
\overline{\Pi}_{V}(w) = F(w)G(w) \quad \quad (w\in {\mathbb{D}}),
\]
for some $G \in L^{\infty}_{B(\cle, \cle_*)}({\mathbb{D}})$. Then
\[
\clh \sim_s \clm\otimes {\mathbb C}^m.
\]
\end{Theorem}

\NI\textsf{Proof.} Note that
\[
ran F (w) = V (w) = ran \overline{\Pi}_{V}(w),
\]
implies that
\[
\overline{\Pi}_{V}(w)F(w) = F(w) \quad \quad (w\in {\mathbb{D}}).
\]
Note also that (cf. (\ref{ineq1}))
\[
\partial
\overline{\Pi}_{V}(w)F(w)=(1-\overline{\Pi}_{V}(w))F^{\prime}(w).
\]
and (cf. Lemma 2.2, \cite{TW})
\[
\partial \overline{\Pi}_{V}(w)=\partial
\overline{\Pi}_{V}(w)\overline{\Pi}_{V}(w),
\]
for all $w \in \mathbb{D}$. We therefore have
$$\begin{array}{lll}
\partial \overline{\Pi}_{V}(w)&=&\partial \overline{\Pi}_{V}(w)\overline{\Pi}_{V}(w)\\
&=&\partial \overline{\Pi}_{V}(w)F(w)G(w)\\
&=&(1-\overline{\Pi}_{V}(w))F^{\prime}(w)G(w)
\end{array}$$
Then for  $M = \|G\|_{\infty} >0$, it follows that
\[
\begin{array}{llll}
\|\partial \overline{\Pi}_{V}(w)\|_2&=&\|(I-\overline{\Pi}_{V}(w))F^{\prime}(w)G(w)\|_2\\
&\leq &\|I-\overline{\Pi}_{V}(w)\|\cdot \|F_2^{\prime}(w)G(w)\|_2\\
&\leq&\|G(w)\|\cdot \|F^{\prime}(w)\|_2\\
&\leq&\|G\|_{\infty}\|F^{\prime}(w)\|_2\\
&= & M\|F^{\prime}(w)\|_2,
\end{array}
\]
for all $w \in \mathbb{D}$. Now as in the proof of Theorem
\ref{similar-link}, we set the bounded subharmonic function $\psi_1$
as $\|M^{\frac{1}{2}}F_2\|^2_2$. Then
\[
trace K_{\clh_1}-trace K_{\clm\otimes {\mathbb{C}}^m}= \|\partial
\overline{\Pi}_{V_1}(w)\|^2_2\leq \Delta \psi_1,
\]
and hence by $(iii) \implies (i)$ of Theorem \ref{similar-link}, we
have that
\[
\clh \sim_s \clm\otimes {\mathbb C}^m.
\]
This completes the proof of the result. \qed

\newsection{Concluding remarks}

A number of questions and directions remain to be explored,
including the similarity problem for the Dirichlet module (however,
see \cite{HK}). We point out that the notion of the Cowen-Douglas
atom does not cover the Dirichlet module (see \cite{St}).

Some of the results of this paper can be generalized in the several
variables set up. However, one of the key ideas to achieve results
of full strength is closely related to the corona problem in several
variables (see \cite{TW}).

Another interesting question relates the quasi-affinity of the
Cowen-Douglas Hilbert modules. For the Hardy space, quasi-affinity
to a submodule of a Hardy module is as same as the Hardy module it
self. It is not known under what additional condition on the frame,
that module will be quasi-affine to a Cowen-Douglas atom of finite
multiplicity.

We also do not known whether the converse of Theorem \ref{last} is
true. More precisely, let $\clh \in B_m(\mathbb{D})$ be a pure
$\clm$-contractive Hilbert module, let $E_{\clh} \cong E_{\clm}
\otimes V$, and let  $V(w) \subseteq \cle$, $w \in \mathbb{D}$, for
some Hilbert space $\cle$. Assume moreover that
$\overline{\Pi}_V$(w) is the orthogonal projection of $\cle$ onto
$V(\bar{w})$. Let
\[
\clh \sim_s \clm \otimes \mathbb{C}^m.
\]
Does there exist a Hilbert space $\cle_*$, a bounded analytic
function $F \in H^{\infty}_{B(\cle_*, \cle)}({\mathbb{D}})$, and a
bounded function $G \in L^{\infty}_{B(\cle, \cle_*)}({\mathbb{D}})$
such that
\[
\mbox{ran~} F(w) = V(w),
\]
and
\[
\overline{\Pi}_{V}(w) = F(w)G(w),
\]
for all $w\in {\mathbb{D}}$?

\vspace{0.2in}

\noindent\textit{Acknowledgement:} The second author is supported in
part by NBHM (National Board of Higher Mathematics, India) grant
NBHM/R.P.64/2014. The second author would also like to thank the
members of the Department of Mathematics of the Hebei Normal
University, China, for their hospitality during a visit during which
part of the work was carried out.

\end{document}